\documentclass[11pt]{amsart} 
\usepackage{latexsym}
\usepackage{amsfonts}
\usepackage{amsmath,amsthm,amssymb}

\newcommand{\cal}{\mathcal}

\newcommand{\CA}{{\cal A}}

\newcommand{\Umin}{U_{\tiny \min}}
\def\epsilon{\varepsilon}
\def\phi{\varphi}

\newcommand{\Aut}{\mbox{Aut}}

\newcommand{\card}{\mbox{card}}

\newcommand{\FN}{F_N}   % F ou F_n ou F_N ?

\newcommand{\FA}{F({\CA})}

\newcommand{\R}{\mathbb R}

\newcommand{\N}{\mathbb N}

\def\strutdepth{\dp\strutbox}
\def \ss{\strut\vadjust{\kern-\strutdepth \sss}}
\def \sss{\vtop to \strutdepth{
\baselineskip\strutdepth\vss\llap{$\diamondsuit\;\;$}\null}}

\def\strutdepth{\dp\strutbox}
\def \sst{\strut\vadjust{\kern-\strutdepth \ssss}}
\def \ssss{\vtop to \strutdepth{
\baselineskip\strutdepth\vss\llap{$\spadesuit\;\;$}\null}}

\def\strutdepth{\dp\strutbox}
\def \ssh{\strut\vadjust{\kern-\strutdepth \sssh}}
\def \sssh{\vtop to \strutdepth{
\baselineskip\strutdepth\vss\llap{$\heartsuit\;\;$}\null}}

\def\tilde{\widetilde}

\def\strutdepth{\dp\strutbox}
\def \ss{\strut\vadjust{\kern-\strutdepth \sss}}
\def \sss{\vtop to \strutdepth{
\baselineskip\strutdepth\vss\llap{$\diamondsuit\;\;$}\null}}

\def\strutdepth{\dp\strutbox}
\def \sst{\strut\vadjust{\kern-\strutdepth \ssss}}
\def \ssss{\vtop to \strutdepth{
\baselineskip\strutdepth\vss\llap{$\spadesuit\;\;$}\null}}
\newtheorem{theorem}{Theorem}%[section]

\theoremstyle{definition}
\newtheorem{definition}[theorem]{Definition}

\newtheorem{remark}[theorem]{Remark}

\theoremstyle{remark}

\numberwithin{equation}{section}
\begin{document}

\author[F.~Ibrahim, M.~Lustig]{Fedaa Ibrahim and Martin Lustig}\address{\tt 
LATP,
Centre de Math\'ematiques et Informatique, 
Aix-Marseille Universit\'e, 
39, rue F.~Joliot Curie, 
13453 Marseille 13, 
France}
\email{\tt fidaa0@hotmail.fr}

\email{\tt Martin.Lustig@univ-amu.fr}

\title[Dual automorphisms of $\FN$]
{Dual automorphisms of free groups}

 
 
\maketitle

{\it To appear in: Extended Conference Abstracts, Spring 2013, CRM Documents, Centre de Recerca Matem\`atica, Bellaterra (Barcelona)}
 
 \bigskip
 
 Throughout this extended abstract of our joint work \cite{FL} we denote by $\FN$ the non-abelian free group of finite rank $N \geq 2$. We assume familiarity of the reader with the basic terminology around the combinatorial theory of automorphisms of free groups.

\smallskip

In order to work with free groups, in particular for algorithmic purposes, one works almost always with a fixed basis $\cal A$ of $\FN$. In this case one has a canonical bijection between elements of $\FN$ and the set $\FA$ of reduced word in $\cal A \cup \cal A^{-1}$.

\smallskip

The distinction between reduced words and free group elements, though unusual in the classical approach to free groups, is rather important for the work presented here. The reason is that a reduced word $w = y_1 \ldots y_r \in \FA$ allows a second, {\em dual} interpretation:  Every such $w$ determines a subset of $\partial \FN$, defined by the {\em cylinder}  $C^1_w \subset \partial \FA$, where by $\partial \FA$ we denote the set of infinite reduced words in $\cal A \cup \cal A^{-1}$:
$$C^1_w = \{x_1 x_2 \dots \mid x_1 = y_1, \ldots, x_r = y_r \}$$
Note that the subset of $\partial \FN$ defined by $C^1_w$ depends not only on the element of $\FN$ given by the word $w$, but also on the chosen basis $\cal A$.  An equivalent way to express this dependency is to note that for any automorphism $\phi \in \Aut(\FN)$ one has in general:
$$\phi(C^1_w) \neq C^1_{\phi(w)}$$
Indeed, the image of a cylinder under an automorphism $\phi \in \Aut(\FN)$ is in general not a cylinder, but a multi-cylinder, i.e. a finite union of cylinders. In his thesis and a subsequent publication (see \cite{Fe}) the first author of the work presented here has given an efficient algorithm how to determine this finite union, and in particular he has proved the following formula:

\begin{theorem}[\cite{Fe}]
\label{Fedaas-formula}
(a) Let $\varphi$ be an automorphism of the free group $\FN$ with finite basis $\cal A$.
For any $u \in F(\cal A)$ there exists a finite set $U \subset F(\cal A)$ such that:
$$\varphi (C^1_u)=C^{1}_U := \bigcup_{u_i \in U} C^1_{u_i}$$

\smallskip
\noindent 
(b) A set $U$ as in statement (a) can be algorithmically derived from $u \in F(\cal A)$ and from the words in the finite subsets $\varphi(\cal A)$ and $\varphi^{-1}(\cal A)$ of $F(\cal A)$. Indeed, the equality in (a) is true for
$$U = \{\varphi(u')|_{S(\varphi)^2} \mid u' \in u|^k\}\, ,$$
with 
$k=S(\varphi)^4+S(\varphi)^3+S(\varphi)^2$, where $S(\varphi)$ is the maximal length of any $\varphi(a_i)$ or $\varphi^{-1}(a_i)$ among all $a_i \in \cal A$.

Here
for any reduced word $w \in F(\cal A)$ and any integer $l \geq 0$ we denote by $w|_l$ the word obtained from $u$ by erasing the last $l$ letters, and by $w|^l$ the set of all reduced words obtained from $u$ by adding $l$ letters from $\cal A \cup \cal A^{-1}$ at the end of $w$.
\end{theorem}

A multi-cylinder $C^1_U \subset \partial F(\cal A)$ does not uniquely define the finite set $U \subset \FA$, but in \cite{Fe} it has been shown that there is a unique {\em reduced subset} $\Umin \subset \FA$ of minimal cardinality which satisfies $C^1_U = C^1_{\Umin}$, and that $\Umin$ can be derived from $U$ by an elementary reduction algorithm. We use these facts to define for any $\phi \in \Aut(\FN)$ the {\em dual automorphism} $\phi^*_\cal A$ as follows:

\begin{definition}
\label{dual-autos}
For any $u \in \FA$ let $\phi^*_\cal A(u)$
be the reduced subset of $\FA$ that satisfies: 
$$C^1_{\phi^*_\cal A(u)} = \phi(C^1_u)$$
\end{definition}

The reader should be warned that, as is indicated by the name, the dual automorphism depends heavily on the choice of the basis $\cal A$ of $\FN$. 

\medskip

There are several immediate natural questions which come to mind if one considers the definition of $\phi^*_\cal A$, concerning for example the seize of the set $\phi^*_\cal A(u) \subset \FN$, its computability, and its behavior under iteration of $\phi^*_\cal A$. We will state now the main results of our joint work \cite{FL}; in particular, Theorem \ref{2n-sets} seems surprising and note-worthy to us.

The following result is inspired by the above Theorem \ref{Fedaas-formula} and proved by similar methods:

\begin{theorem}
\label{finiteness}
For any automorphism $\phi \in \Aut(\FN)$ and any basis $\cal A$ of $\FN$ there exists a finite collection $\cal U(\phi^*_\cal A) = \{U_1, \ldots, U_r \}$ of finite sets $U_i \subset \FA$ of reduced words, such that for every $w \in \FA$ the dual image $\phi^*_\cal A(w)$ is given by $v U_i$, for some $v \in \FA$ and $U_i = U_i(w) \in \cal U(\phi^*_\cal A)$. The word $v$ can be specified further to $v = \phi(w|_k)|_K$ for some constants $k, K \geq 0$ which are independent of $w$.

In particular, the cardinality of $\phi^*_\cal A(w)$, for any $w \in \FA$, is bounded above by a constant which only depends on $\phi$ (and the fixed basis $\cal A$).
\end{theorem}

The computation of the constants $k$ and $K$ and of the finite sets $U_i$ is possible by the use of Theorem \ref{Fedaas-formula}, and writing an efficient computer program should not be a very difficult task.

\smallskip

However, the result stated below in Theorem \ref{2n-sets} is much more striking and also more useful for computational purposes; it came about when we tried to look at the easiest non-trivial example, namely an elementary Nielsen automorphism:

\begin{remark}
\label{fundamental-formula}
[Fundamental formulas]
\rm
\, We recall Nielsen's famous result that every automorphism of a free group can be written as product of {\em elementary automorphisms}, which can be regrouped into the following two types:

\smallskip
\noindent
(a)
For the ``standard'' the elementary Nielsen automorphism $\phi \in F(\cal A)$,  for $\cal A = \{a, b, c_1, \ldots, c_q\}$, given by $a \mapsto ab, \, b \mapsto b, c_j \mapsto c_j$, we obtain (where we use the convention that a product of words $u v$ is written as $u \cdot v$ if no cancellation occurs between the end of $u$ and the beginning of $v$):
$$\phi^*_\cal A(w \cdot a) = \{\phi(w) \cdot a\}$$
$$\phi^*_\cal A(w \cdot a^{-1}) = \{\phi(w) b^{-1} \cdot a^{-1})\}$$
$$\phi^*_\cal A(w \cdot b) = \{\phi(w)\cdot b, \, \phi(w)\cdot a^{-1}\}$$
$$\phi^*_\cal A(w \cdot b^{-1}) = \{\phi(w) b^{-1} \cdot b^{-1}, \phi(w) b^{-1} \cdot a\}$$
$$\phi^*_\cal A(w \cdot c_j) = \{\phi(w) \cdot c_j\}$$
$$\phi^*_\cal A(w \cdot c_j^{-1}) = \{\phi(w) \cdot c_j^{-1}\}$$

\smallskip
\noindent
(b)
The other class of elementary automorphisms $\phi$ is given by a permutation of the letters in $\cal A$, or by an inversion of some of them.
In other words, one has $|\phi(a_i)| = 1$ for every $a_i \in \cal A$, and there is never a cancellation in the image of a reduced word.
As a consequence one obtains directly
$$\phi^*_\cal A(w) = \{\phi(w)\}$$
for any $w \in F(\cal A)$.
\end{remark}

These fundamental formulas give a straight forward method to calculate the analogous formulas for any automorphism, by applying the above fundamental formulas in an iterative way.  
This leads directly to:

\begin{theorem}
\label{2n-sets}
For any automorphism $\phi \in \Aut(\FA)$ the collection $\cal U(\phi^*_\cal A)$ from Theorem \ref{finiteness} can be specified to consist of precisely $2N$ sets $U(x)$, with $x \in \cal A \cup \cal A^{-1}$, such that for every reduced word $w = y_1 \ldots y_q$ one has (using the terminology of Theorem \ref{finiteness}):
$$U_i(w) = U(y_q)$$
Furthermore, if $\phi$ is the product of elementary automorphisms (i.e. basis permutations, basis inversions, or elementary Nielsen automorphisms), then the cardinality of each $U(x)$ is bounded by $2^t$, where $t$ is the number of elementary Nielsen automorphism in the above decomposition of $\phi$, and the constants $k \geq 0$ and $K \geq 0$ from Theorem \ref{finiteness} can be chosen as $k = 1$ and $K = 0$.
\end{theorem}

The $2^t$-bound from the last theorem suggests the definition of a {\em dual growth rate} $\lambda^*_{\phi}$ which can be defined for example as the limit superior of the family of 
values $(\card (\phi^k)^*_\cal A(x))^{\frac{1}{k}}$,
for $x \in \cal A \cup \cal A^{-1}$ and $k \in \N$.

\begin{theorem}
\label{dual-growth-rate}
For any automorphism $\phi \in \Aut(\FN)$ the dual growth rate $\lambda^*_{\phi}$ is independent of the choice of the basis $\cal A$. It can be calculated as Perron-Frobenius eigenvalue of a $2N \times 2N$-matrix which can be algorithmically derived from $\phi$.
\end{theorem}

\medskip

It turns out that cylinders are less natural objects in group theory than in combinatorics: For example, on has $vC^1_u = C^1_{vu}$ only if $u$ isn't completely cancelled when reducing $vu$. This (and related) problems vanish if one passes over to {\em double cylinders} $C^2_{[v, w]}$, for $v, w \in \FA$, which are defined as sets of endpoint pairs of biinfinite reduced paths in the Cayley tree $\tilde \Gamma_\cal A(\FN)$ that pass through both vertices $v$ and $w$. All of the result stated above for single cylinders have natural analogues for double cylinders (compare \cite{Fe}). Double cylinders play a natural role in defining currents for free groups and in properly setting up the basic theory of such, see \cite{K1, K2}. Indeed, one way to name concretely a current $\mu$ is to specify its Kolmogorov function, which is the non-negative function $\mu_\cal A: \FA \to \R$ which associates to every $w \in \FA$ the current measure $\mu(C^2_{[1, w]})$ of the double cylinder $C^2_{[1, w]} \subset \FN \times \FN \smallsetminus {\rm diagonal}$.

\smallskip

In \S 6 of \cite{K2} a formula is derived how to derive from the Kolmogorov function for a current $\mu$ that for the image current $\phi(\mu)$, for any $\phi \in \Aut(\FN)$.  The above finiteness statements and the algorithmic applications allow a substantial simplification and specification of this formula. The interest towards a concrete and effective calculation of $\phi(\mu)$ from known data for $\phi$ and $\mu$ was one of the original motivations to stimulate this work.

\end{document}